\documentclass[11pt]{amsart}

\usepackage[utf8]{inputenc}
\usepackage[english]{babel}
\usepackage{amsmath, amssymb, amsthm}
\usepackage{mathtools}
\usepackage{tikz-cd}
\usepackage{xcolor}
\usepackage[
    colorlinks=true,
    linkcolor={blue!60!black},
    citecolor={blue!60!black},
    urlcolor={blue!80!black}
]{hyperref}
\usepackage[capitalize]{cleveref}
\usepackage{enumitem}
\usepackage{float}

\theoremstyle{plain}
\newtheorem{theorem}{Theorem}[section]
\newtheorem{lemma}[theorem]{Lemma}
\newtheorem{proposition}[theorem]{Proposition}

\theoremstyle{definition}
\newtheorem{definition}[theorem]{Definition}

\newtheorem{question}[theorem]{Question}
\theoremstyle{remark}
\newtheorem{remark}[theorem]{Remark}
\theoremstyle{convention}
\newtheorem{convention}[theorem]{Convention}

\newcommand{\N}{\mathbb{N}}   
\newcommand{\Z}{\mathbb{Z}}   
\newcommand{\Q}{\mathbb{Q}}

\newcommand{\G}{\mathcal{G}}
\newcommand{\x}[1]{x_{#1}(-)}
\newcommand{\hr}{\theta}

\title[Steinberg module and top-dimensional cohomology]{Steinberg module and top-dimensional cohomology for arithmetic Chevalley groups}
\author{Miguel Gentili}
\address{Department for Mathematical Logic and Foundational Research, Münster University, Germany}
\email{miguel.gentili@uni-muenster.de}
\date{}

\begin{document}
\maketitle

\begin{abstract}
We prove that for a Chevalley--Demazure group scheme $\G$ of type $G_2$, $F_4$ or $E_8$ and a Euclidean ring $\mathcal{O}$ with field of fractions $K$, the Steinberg module of $\G(K)$ is generated by integral apartment classes. As a consequence, the top-dimensional rational cohomology of $\G(\mathcal{O})$ vanishes if $\mathcal{O}$ is a ring of integers, extending results of Tóth and of Brück--Santos Rego--Sroka.
\end{abstract}

\section{Introduction}
In \cite{CFP14}, Church, Farb and Putman conjectured that the rational cohomology $H^{\operatorname{vcd}(\mathrm{SL}_n(\Z))-i}(\mathrm{SL}_n(\Z); \Q)$ of the special linear group over the integers vanishes for $i < \operatorname{rk}(\mathrm{SL}_n)=n-1$. A proof has recently been announced in the preprint~\cite{Den26}. Motivated by this conjecture, Brück, Santos Rego and Sroka asked in \cite{BSS24} whether a similar vanishing pattern holds for other Chevalley groups.

\begin{question}[Brück–Santos Rego–Sroka {\cite[Question 1.2]{BSS24}}]
Let $\G$ be a semisimple Chevalley--Demazure group scheme and let $\mathcal{O}$ be the ring of integers of a number field. If $\mathcal{O}$ is Euclidean, is it true that 
\[
H^{\operatorname{vcd}(\G(\mathcal{O}))-i}(\G(\mathcal{O}),\Q)=0
\]
for all $i < \operatorname{rk}(\G)$?
\end{question}

In the same article, the authors obtained positive answers to this question for $i=0$ and for Chevalley--Demazure group schemes of type $A_n$, $B_n$, $C_n$, $D_n$, $E_6$ and $E_7$. Previous results were known for types $A_n$ and $C_n$. Namely, Lee and Szczarba \cite{LS76} proved it for type $A_n$ and Brück, Patzt and Sroka \cite[Chapter 5]{Sro21} for type $C_n$.

In this article we complete the case $i=0$ by treating the remaining types $G_2$, $F_4$ and $E_8$. More precisely, we prove the following theorem:

\begin{theorem} \label{the:theorem}
Let $\mathcal{O}$ be the ring of integers of a number field, and let $\G$ be a semisimple Chevalley--Demazure group scheme with irreducible root system. If $\mathcal{O}$ is Euclidean, then
\[
H^{\operatorname{vcd}(\G(\mathcal{O}))}(\G(\mathcal{O}),\Q)=0.
\]
\end{theorem}

The proof follows the approach used by Brück, Santos Rego and Sroka. Their argument relies on two main ingredients: Borel--Serre duality \cite{BS73} and a paper of Tóth \cite{Tot05}. Throughout, $\mathcal{O}$ denotes the ring of integers of a number field. Whenever the argument does not use this assumption, we write $O$ for an arbitrary integral domain instead.

If $\G$ is a Chevalley--Demazure group scheme and $K$ is a field, to $\G(K)$ one associates the so-called \emph{Steinberg module} $\mathrm{St}$, on which $\G(K)$ acts, together with a distinguished element $[\Sigma]$ called the \emph{standard apartment class}. The Solomon--Tits \Cref{solomon:tits:thm} states that the Steinberg module is generated by $\G(K)$-translates of $[\Sigma]$, i.e. by \emph{apartment classes}. If $O$ is an integral domain with field of fractions $K$, we call \emph{integral apartment classes} the $\G(O)$-translates of $[\Sigma]$. Brück, Santos Rego and Sroka proved the following:

\begin{theorem}[{Brück--Santos Rego--Sroka, \cite[Theorem 3.4]{BSS24}}]  \label{vanishing:cohomology}
Let $K$ be a number field with ring of integers $\mathcal{O}$ and let $\G$ be a semisimple Chevalley--Demazure group with irreducible root system. If the Steinberg module of $\G(K)$ is generated by integral apartment classes, then
\[
H^{\operatorname{vcd}(\G(\mathcal{O}))}(\G(\mathcal{O}); \Q) = 0.
\]
\end{theorem}

The proof of this result relies on Borel--Serre Duality \ref{BorelSerre}. Borel and Serre proved that arithmetic groups, such as $\G(\mathcal{O})$, satisfy a duality between rational cohomology and homology with coefficients in the $\Q[\G(\mathcal{O})]$-module $\Q \otimes_\Z \mathrm{St}$. Therefore, the problem of computing top--dimensional cohomology can be reduced to computing co--invariants of $\Q \otimes_\Z \mathrm{St}$.

In the aforementioned paper \cite{Tot05}, Tóth proved that for types $A_n$, $B_n$, $C_n$, $D_n$, $E_6$, and $E_7$, the Steinberg module is indeed generated by integral apartment classes if the ring $O$ is Euclidean. In these cases, the existence of a maximal parabolic subgroup whose unipotent radical is abelian simplifies some computations. The remaining types admit no such parabolic subgroup, and it is for this reason that they were excluded by Tóth. Our main contribution in this article is the following theorem, which shows that the same techniques nonetheless apply to these types.

\begin{theorem} \label{integral:generation}
Let $O$ be a Euclidean domain with fraction field $K$, and let $\G$ be a semisimple Chevalley--Demazure group scheme associated to a root system of type $G_2$, $F_4$ or $E_8$. Then, the Steinberg module for $\G(K)$ is generated by integral apartment classes.
\end{theorem}

Combining \Cref{vanishing:cohomology} with the theorem above yields \Cref{the:theorem} for types $G_2$, $F_4$ and $E_8$; the remaining types were established in \cite{BSS24}. This completes the proof of \Cref{the:theorem}.

We split the proof of \Cref{integral:generation} into two separate sections: in \Cref{local:section} we prove a ``local'' version of this fact (\Cref{local:thm}), which will be necessary for the general proof contained in \Cref{global:section}.

We point out that the assumption that $\mathcal{O}$ is Euclidean cannot be dropped entirely. For types $A_n$ and $C_n$ the weaker assumption that $\mathcal{O}$ is a PID is necessary: it was proved that the cohomology does not vanish when $\mathcal{O}$ is not a PID (see \cite[Theorem D]{CFP19} and \cite[Theorem 1.1]{BH23} for types $A_n$ and $C_n$ respectively).

\subsection*{Acknowledgements}
I would like to thank Benjamin Br\"uck for suggesting this problem to me
and for many helpful discussions.
This work is funded by the Deutsche Forschungsgemeinschaft (DFG, German Research Foundation) -- Project-ID 561029452, and under Germany's Excellence Strategy EXC 2044/2 -- 390685587, Mathematics M\"unster: Dynamics--Geometry--Structure.

\section{Background}

\subsection{Chevalley groups} \label{sec:chevgps}
In this section we briefly review the theory of \emph{Chevalley groups}. This class of algebraic groups was first introduced by Chevalley in his seminal paper \cite{Che55} and later studied by Ree \cite{Ree64} and Kostant \cite{Kos66}.
A comprehensive treatment of such groups via modern algebraic
geometry was given by Demazure and Grothendieck in \cite{DG70}. We refer to it for the theory over an arbitrary base ring, and to the
more accessible \cite{Mil17} for the theory over a field. Precisely, we have the following.

\begin{definition}
A \emph{Chevalley--Demazure group scheme} is a split reductive affine
group scheme $\G$ over $\Z$. For a commutative ring $R$ with unit, the
associated \emph{Chevalley group} is the group $\G(R)$ of its
$R$-points.
\end{definition}

That is, a Chevalley--Demazure group scheme is a representable functor from the category of rings to the category of groups satisfying some additional properties. We will restrict to semisimple Chevalley--Demazure group schemes associated to $G_2$, $F_4$ and $E_8$.
In some situations, it suffices to think of a Chevalley group as an
abstract group given by explicit generators and relations. 
Let $R$ be a ring and consider the elements of the Chevalley group $\G(R)$ given by the symbols $\{x_{\alpha}(t) \mid \alpha \in \Phi,\ t \in R\}$, where $\Phi$ denotes the root system associated to $\G$. These elements satisfy the following relations:
\begin{enumerate}[label=\textbf{(R\arabic*)}, ref={[R\arabic*]}]
    \item\label{rel:additive}
    $x_{\alpha}(t_1)\, x_{\alpha}(t_2) = x_{\alpha}(t_1 + t_2)$,
    for all $\alpha \in \Phi$ and $t_1, t_2 \in R$;

    \item\label{rel:chevform}
    (\textbf{Chevalley Commutator Formula}) for all linearly
    independent $\alpha, \beta \in \Phi$ and $t, u \in R$,
    \[
        [x_{\alpha}(t),\, x_{\beta}(u)]
        = \prod_{\substack{i,j > 0 \\ i\alpha + j\beta \in \Phi}}
          x_{i\alpha + j\beta}\!\left(C_{ij}^{\,\alpha\beta}\, t^{i} u^{j}\right),
    \]
    where the product is taken over the roots $i\alpha + j\beta$ in a fixed order, and the $C_{ij}^{\,\alpha\beta}$ are certain integers depending on the root system $\Phi$;

    \item\label{rel:weyl}
    $w_{\alpha}(t)\, x_{\alpha}(u)\, w_{\alpha}(t)^{-1}
    = x_{-\alpha}\!\left(-t^{-2} u\right)$,
    for all $t \in R^{\times}$ and $u \in R$, where
    $w_{\alpha}(t) \coloneqq x_{\alpha}(t)\, x_{-\alpha}(-t^{-1})\, x_{\alpha}(t)$;

    \item\label{rel:torus}
    $h_{\alpha}(t)\, h_{\alpha}(u) = h_{\alpha}(tu)$,
    for all $t, u \in R^{\times}$, where
    $h_{\alpha}(t) \coloneqq w_{\alpha}(t)\, w_{\alpha}(-1)$.
\end{enumerate}

In the Chevalley Commutator Formula \ref{rel:chevform}, we order the terms in increasing $i+j$.  Note that if $i_1+j_1=i_2+j_2$, then $x_{i_1\alpha+j_1 \beta}$ and $x_{i_2 \alpha+j_2 \beta}$ commute.

If $R$ is a Euclidean domain and $\G$ is simply connected, the relations above no longer form a complete set, but $\G(R)$ is still generated by the root elements $x_\alpha(t)$ with $\alpha\in\Phi$ and $t\in R$
\cite[Corollary~3 to Theorem~18, Chapter~8]{Ste16}. This may fail if $R$ is
not Euclidean; see \cite{Coh66} for an example.

The following facts about split reductive groups can be found in \cite[Exposé~XXII]{DG70}.
For every $\alpha \in \Phi$, there is a closed embedding
$U_\alpha \hookrightarrow \G$ of the \emph{root subgroup} $U_\alpha$,
whose $R$-points are given by
$U_\alpha(R) \coloneqq \langle x_\alpha(t) \mid t \in R \rangle$.
The map $t \mapsto x_\alpha(t)$ induces an
isomorphism of group schemes $U_\alpha \cong \mathbb{G}_a$, where
$\mathbb{G}_a$ is the additive group scheme, defined by
$\mathbb{G}_a(R) \coloneqq (R,+)$ for every commutative ring $R$ with unit.

By definition, every Chevalley–Demazure group scheme $\G$ contains a split maximal torus $T$. The \emph{Weyl group} of $\G$ is the quotient group scheme $W \coloneqq N_{\G}(T)/T,$ where $N_{\G}(T)$ denotes the normalizer of $T$ in $\G$. It is a finite constant group scheme isomorphic to the abstract Weyl group $W_\Phi = \langle s_\alpha \mid \alpha \in \Delta \rangle$ of the root system $\Phi$, generated by the reflections $s_\alpha$ associated to the simple roots $\alpha \in \Delta$ (cf.\ \cite[Proposition~21.1, Theorem~21.37]{Mil17}, \cite[Exposé~XXII.3]{DG70}). 

Let $R$ be a commutative ring with unit. For each simple root $\alpha \in \Delta$, set
\[
  w_\alpha \coloneqq x_\alpha(1)\,x_{-\alpha}(1)^{-1}\,x_\alpha(1) \in N_{\G}(T)(R) \subseteq \G(R),
\]
a distinguished representative of the preimage of $s_\alpha$ under the isomorphism above. We call $\tilde{W} \coloneqq \langle w_\alpha \mid \alpha \in \Delta \rangle \subseteq \G(R)$ the \emph{extended Weyl group} of $\Phi$.

\subsection{Tits buildings and Steinberg modules}

We refer to Abramenko and Brown's book \cite{AB08} for a comprehensive introduction to the theory of Tits buildings. Let $G$ be an algebraic group over a field $K$ with a $BN$-pair; the subgroup $B$ is called \emph{standard Borel subgroup of G}.  For instance, let $G=\G$ be a Chevalley--Demazure group scheme associated to a root system $\Phi$ with basis $\Delta$. Then $\G$ has a $BN$-pair for any field $K$: with $T$ the fixed split maximal torus and $U_\alpha$ the unipotent subgroups of \Cref{sec:chevgps}, one takes $B = \langle T(K), U_\alpha(K) \mid \alpha \in \Phi^+ \rangle$ and $N = N_{\G}(T)(K)$ the normalizer of $T$, where $\Phi^+$ is the set of positive roots.

Any subgroup $P$ containing the standard Borel $B$ is called \emph{standard parabolic}. The \emph{spherical building} $\Delta(\G,K)$ associated to $\G(K)$ is a simplicial complex whose vertices are right cosets of maximal proper standard parabolic subgroups, and $k+1$ vertices $ P_1 g_1, \dots, P_{k+1} g_{k+1}$ form a $k$-simplex if their intersection is non--empty (i.e.\ they all contain a fixed $B g$ for some $g \in \G(K)$). Note that $\G(K)$ acts on $\Delta(\G,K)$ by right translation, and hence so does $\G(R)$ for every subring $R\subseteq K$.

The maximal simplices are called \emph{chambers}, and the standard Borel $B$ is called the \emph{fundamental chamber}. The \emph{standard apartment} is the subcomplex $\Sigma \coloneqq \langle B \tilde{w} \mid \tilde{w} \in \tilde{W} \rangle$, where $\tilde{W}$ is the extended Weyl group of $\Phi$. 

There exists a distinguished class in homology, namely
\[
[\Sigma] \coloneqq \sum_{w \in W_\Phi} (-1)^{\ell(w)} \{B \tilde{w}\} \in \tilde{H}_{\operatorname{rk}(\Phi)-1}(\Delta(\G,K);\Z),
\]
where $\ell(w)$ is the length of $w$ and $\tilde{w}$ is a representative of $w \in W_\Phi$ in the extended Weyl group $\tilde{W}$.

\begin{definition}
Let $\G$ be a Chevalley--Demazure group scheme over a field $K$, and let $\Delta(\G,K)$ be the associated Tits building. The \emph{Steinberg module} of $\G(K)$ is the $\Z[\G(K)]$-module defined as 
\[
\mathrm{St}(\G, K) = \mathrm{St}\coloneqq \tilde{H}_{\operatorname{rk}(\Phi)-1}(\Delta(\G,K);\Z).
\]
An element of the form $[\Sigma]g$ for $g \in \G(K)$ is called an \emph{apartment class}. If $K$ is the field of fractions of an integral domain $O$, the apartment class $[\Sigma]g$ is called \emph{integral} if $g \in \G(O)$.
\end{definition}

\begin{theorem}[{Solomon--Tits, \cite{Sol69}[Theorem~4.73], \cite{AB08}}] \label{solomon:tits:thm}
The Steinberg module is a free abelian group. Apartment classes form a generating set, but not a basis. In particular, the map of right $\Z[\G(K)]$-modules
\[
\pi : \Z[\G(K)] \rightarrow \mathrm{St},\qquad g \mapsto [g]=[\Sigma]\cdot g,
\]
is surjective, where $[\Sigma]$ is the fundamental class of the standard apartment.
\end{theorem}

The Solomon--Tits Theorem says that, as a right $\Z[\G(K)]$-module, the Steinberg module is isomorphic to the quotient $\Z[\G(K)] / \ker (\pi)$, where $\pi$ is the map above. It is interesting to understand what $\ker(\pi)$ is exactly.

\begin{definition} \label{def:R}
Let $e$ be the unit element of $\G(K)$. We denote by $\mathcal{R}$ the right $\Z[\G(K)]$-submodule generated by
\[
\{e-x_\alpha(t) - x_{-\alpha}(t^{-1}), e-h, e+w_\alpha \mid \alpha \in \Phi, t \in K^*, h \in T(K) \}.
\]
\end{definition}

\begin{proposition}[{Tóth, \cite[Theorem 2]{Tot05}}] \label{R:ker}
Let $\pi$ be the projection map of Solomon--Tits \Cref{solomon:tits:thm} and let $\mathcal{R}$ be the module defined in \Cref{def:R}. Then
\[
\mathcal{R} = \ker(\pi).
\]
\end{proposition}

\subsection{Arithmetic groups}
We briefly recall the notion of \emph{virtual cohomological dimension}. More details can be found in \cite[Chapter~VIII]{Bro82}.

Given a group $\Gamma$, its \emph{cohomological dimension} is
\[
\operatorname{cd}(\Gamma) \coloneqq \sup \{ n \mid H^n(\Gamma;M) \neq 0 \text{ for some $\Gamma$-module } M \},
\]
and we set $\operatorname{cd}(\Gamma) \coloneqq \infty$ if the supremum is not finite. For instance, if $\Gamma$ has torsion, then $\operatorname{cd}(\Gamma) = \infty$ (see \cite[Chapter~VIII, Corollary~2.5]{Bro82}).

The \emph{virtual cohomological dimension} of a group $\Gamma$ is the cohomological dimension of any torsion-free subgroup of finite index of $\Gamma$. By a theorem of Serre \cite{Ser71}, any two such subgroups have the same cohomological dimension, so this is a well-defined invariant of $\Gamma$. Moreover, by Selberg's lemma, every arithmetic group contains a torsion-free subgroup of finite index (cf.\ \cite{Sel60}, \cite[Corollary~7.6.4]{Rat19}).

Let $K$ be a number field with ring of integers $\mathcal{O}$. If $\G$ is a Chevalley--Demazure group, then $\G(\mathcal{O})$ is an arithmetic group (see \cite[Chapter~I.3]{Mar91}).
The group $\G(\mathcal{O})$ acts properly and cocompactly on the Borel--Serre bordification of the associated symmetric space, whose boundary is homotopy equivalent to the spherical building $\Delta(\G,K)$. Exploiting this action, Borel and Serre proved that arithmetic groups are \emph{virtual duality groups} in the sense of Bieri--Eckmann \cite{BE73}; that is, every torsion-free finite-index subgroup of an arithmetic group satisfies a generalized Poincar\'e--Lefschetz duality, with dualizing module the Steinberg module $\mathrm{St}$. More precisely, we have the following:
\begin{theorem}[Borel--Serre Duality, \cite{BS73}] \label{BorelSerre}
Let $\mathcal{O}$ be the ring of integers of a number field, and let $\G$ be a Chevalley--Demazure group. Then the virtual cohomological dimension of $\G(\mathcal{O})$ is finite. Moreover, for every $i$, there is an isomorphism
\[
H^{\operatorname{vcd}(\G(\mathcal{O}))-i}(\G(\mathcal{O}); \Q) \;\cong\; H_i\!\left(\G(\mathcal{O});\, \Q \otimes_{\Z} \mathrm{St}\right),
\]
where $\mathrm{St}$ denotes the Steinberg module.
\end{theorem}
We need not pass to torsion-free subgroups, since we use rational coefficients.

\section{Proof of \Cref{integral:generation}}

\subsection{The local statement} \label{local:section}
In this section, $\Phi$ denotes a root system with basis $\Delta$, while $\G$ is a semisimple Chevalley--Demazure group scheme and $U$ is the unipotent radical of its Borel. The ring $O$ is an integral domain with field of fractions $K$. The group of units of a ring $R$ is denoted by $R^{\times}$. We denote by $\mathcal{R}$ the relations module of \Cref{def:R}. Many of the proofs and the arguments are taken from \cite{Tot05}.

Let $\bar{r} \in \Delta$ be a distinguished simple root. The $\bar{r}$-height of a root $\alpha \in \Phi$ is the coefficient of $\bar{r}$ in the unique expression of $\alpha$ as a combination of simple roots. The groups considered by Tóth contain maximal parabolic subgroups whose unipotent radical is abelian, which is equivalent to the existence of a simple root $\bar{r}$ such that the $\bar{r}$-height of any root is $-1$, $0$, or $1$ (\cite[Lemma 2.4]{RRS92}). This simplifies some calculations and is the reason why $G_2$, $F_4$, and $E_8$ are excluded in Tóth's argument. The purpose of this article is to extend Tóth's result to these groups. Despite some technical difficulties that arise, we will do so by following his approach.

From now on, let $\Phi \in \{G_2, F_4, E_8 \}$ with basis $\Delta$. Looking at the plates in Bourbaki's book \cite[Plates VII--IX]{Bou02}, there exists a root $\bar{r} \in \Delta$ such that the $\bar{r}$-height of every positive root is less than or equal to $2$. Moreover, the maximal root is the only root of $\bar{r}$-height $2$. Namely, for $G_2$ we call the short and the long simple root $\alpha$ and $\beta$ respectively, and we take $\bar{r} = \beta$. For $E_8$ we take $\bar{r}= \alpha_8$ and for $F_4$ we pick $\bar{r}=\alpha_1$, where simple roots are labelled as with Bourbaki's convention. From now on, we will refer to the $\bar{r}$-height of a root simply as its \emph{height}. Note that the word ``height'' has a different standard meaning in this context, namely the sum of the coefficients appearing in the expression of a root as a combination of simple roots. We will use this notion in \Cref{root:order} to define an order on the lattice spanned by the root system, and never mention it elsewhere, so no confusion should arise. We will denote by $\Phi_1$ the roots of height $1$, and by $\Phi_{>0}$ the roots of positive height. For example, if $\Phi = G_2$, we have $\Phi_1 = \{\beta, \alpha + \beta, 2\alpha + \beta, 3\alpha + \beta \}$ and $\Phi_{>0} = \Phi_1 \cup \{3\alpha + 2 \beta\}$.

\begin{definition}
Let $O$ be an integral domain. A prime $p \in O$ is called \emph{small} if the map
\[
O^\times \to (O/(p)) \setminus \{0\}
\]
is surjective.
\end{definition}

For instance, if $O = \Z$, the prime $5$ is not small: indeed, $[2] \in \Z/(5)$ cannot be represented by a unit. The only positive small primes in $\Z$ are $2$ and $3$.

Our local statement is the following proposition. It is ``local'' in the sense that it will allow us to treat one prime at a time in the proof of \Cref{integral:generation}.
\begin{proposition} \label{local:thm}
Let $\G$ be a semisimple Chevalley--Demazure group scheme, let $O$ be an integral domain with fraction field $K$, let $p$ be a small prime and let $O' \coloneqq O[p^{-1}]$. Then 
\[
\Z[U(O ')] \subseteq  \Z[\G(O)] + \mathcal{R}.
\]
\end{proposition}

All the cases except types $G_2$, $F_4$ and $E_8$ are treated in \cite[Theorem 3]{Tot05}. In the rest of this section we complete the remaining cases.

Fix the distinguished simple root $\bar{r} \in \Delta$ for which the height of any root is less than or equal to $2$. The group
\[
\G' \coloneqq \langle U_s, U_{-s} \rangle_{s \in \Delta \setminus \{\bar{r}\}}
\]
is a Chevalley--Demazure group scheme of rank $\operatorname{rk}(\G) - 1$, generated by the root elements associated to roots of height $0$. For $\G$ of type $G_2$, $F_4$ and $E_8$, this corresponds to $A_1$, $C_3$ and $E_7$ respectively, and the result holds for these types.

For $R$ either $O$ or $O'$, let $A(R)$ and $B(R)$ be the sets of products in a fixed order of root elements of height $1$ and of height $2$, respectively. 
Since $\G'(R)$ is generated by root elements of height $0$, it normalizes $A(R)B(R)$ by the Chevalley Commutator Formula.

In view of the next lemma, we can reduce to proving the statement for elements of $A(O')B(O')$. The proof follows Tóth's argument in \cite[Lemma 2]{Tot05} with minor modifications; we reproduce it for completeness.
\begin{lemma} \label{height 0}
Let $\G$ be a group of type $G_2$, $F_4$ or $E_8$. Then
\[
\Z[U(O')] \subseteq \langle \Z[A(O')B(O')] \rangle_{\Z[\G(O)]} + \mathcal{R},
\]
where $\langle \Z[A(O')B(O')] \rangle_{\Z[\G(O)]}$ is the right $\Z[\G(O)]$-module generated by $\Z[A(O')B(O')]$.
\end{lemma}
\begin{proof}
Rearranging root elements using the Chevalley Commutator Formula, we can express $u \in U(O')$ as
\[
u=\lambda ab,
\]
where $\lambda \in \G'(O') \cap U(O')$, $a \in A(O')$, and $b \in B(O')$. Since $\lambda$ lies in a Chevalley group for which \Cref{local:thm} holds, we can write $\lambda = \sum_i n_i \lambda_i + \rho$, where $\lambda_i \in \G'(O)$, $n_i \in \Z$ and $\rho \in \mathcal{R'}$, where $\mathcal{R'}$ is the relation module of \Cref{def:R} for $\G'$. This decomposition still holds if we view the $\lambda_i$ as elements of $\G(O)$. We have
\[
u = \sum_i n_i\,(\lambda_i a  b \lambda_i^{-1})\cdot \lambda_i + \rho\, ab;
\]
since $\G'(O')$ normalizes $A(O')B(O')$, we have $\lambda_i a  b \lambda_i^{-1} \in A(O')B(O')$. Moreover, $\mathcal{R}' \subseteq \mathcal{R}$, and the claim follows.
\end{proof}

We introduce a valuation that, roughly speaking, keeps track of how far an element in $O'$ is from lying in $O$. Namely, for $k \in O$, define
\[
\nu_p(k) \coloneqq \sup\{n \geq 0 : k \in (p^n) \},
\]
where the supremum is taken to be $\infty$ if $k \in (p^n)$ for every $n \in \N$. For $t = k/p^n \in O'$, define $\nu_p(t) \coloneqq \nu_p(k) - n$. Note that $t \in O$ if and only if $\nu_p(t) \geq 0$. We denote by the same symbol the map
\[
\nu_p \colon A(O')B(O') \longrightarrow (\Z \cup \{\infty\}) \times \Phi_{>0}
\]
defined as follows. Given $ab = \prod_r x_r(t_r)$, let $n \coloneqq \min_r \nu_p(t_r)$ be the smallest valuation among the $t_r$, and let $s$ be the smallest root for which $\nu_p(t_s) = n$. We then set $\nu_p(ab) \coloneqq (n, s)$. We order the so-obtained tuples by lexicographic order. Note that to order the roots we need to fix a total order on $\Phi_{>0}$. To do so, we adopt the following:

\begin{convention} \label{root:order}
We fix a total order on the root lattice spanned by the root system $\Phi$ and which is compatible with addition and for which every positive root is greater than $0$.
For instance, pick the order of simple roots as in Bourbaki. Every element of the root lattice can be uniquely written as combination of simple roots. One may compare two elements of the root lattice first by the sum of the coefficients in their expression as a combination of simple roots, and break ties lexicographically on the coefficient vector. The set
$\Phi_{>0}$ is ordered by restriction, and it is with respect to this
order that the second entry of $\nu_p$ is computed.  

For example, if $\alpha_1,\alpha_2$ are simple roots labelled as in
Bourbaki, then $3\alpha_1+\alpha_2 > \alpha_1+3\alpha_2$: the two elements
have the same sum of coefficients, and in the tie-break the coefficient of
$\alpha_1$ is compared first.
\end{convention}

Throughout, we denote with $\hr$ the unique root of height $2$. Let $ab = \prod_{r\in \Phi_1} x_r(t_r) x_\hr(t_\hr)$ be an element of $A(O')B(O')$, where $a = \prod_{r \in \Phi_1} x_r(t_r) \in A(O')$ and $b = x_\hr(t_\hr) \in B(O')$. The idea of the proof is to replace $ab$ by an element of higher valuation, by suitably adding an element of $\mathcal{R}$, until the valuation is at least $0$, so that the resulting element belongs to $\Z[\G(O)] + \mathcal{R}$.

\begin{remark}
A remark on the notation that we will use: a prime ($'$) indicates that an element has been modified in a way that does not decrease its valuation; depending on the context, the valuation is either strictly increased or merely preserved. Moreover, across different lemmas we reuse the same letter for different elements; this causes no confusion, since within each statement only the valuation of the element is relevant.
\end{remark}

\begin{remark} \label{swapping:elements}
Before proceeding, we need to investigate more carefully what happens when we swap two root elements. By the Chevalley Commutator Formula, elements of height $2$ commute with each other and with every element of height $1$. On the other hand, the commutator of elements of height $1$ is, in general, not trivial. Let $r_1, r_2 \in \Phi_1$. By the Chevalley Commutator Formula, their commutator has at most one factor, namely $x_{r_1 + r_2}\!\left(C^{r_1, r_2}_{1,1} t_{r_1} t_{r_2}\right)$. Since the constants appearing in the formula are integers, we have
\[
\nu_p\!\left(C^{r_1, r_2}_{1,1} t_{r_1} t_{r_2}\right) \geq \nu_p(t_{r_1} t_{r_2}) = \nu_p(t_{r_1}) + \nu_p(t_{r_2}) \geq 2 \min\{\nu_p(t_{r_1}), \nu_p(t_{r_2})\}.
\]
That is, swapping two root elements of height $1$ produces a root element of height $2$ whose valuation is at least twice the smaller of the two valuations. Moreover, since $\nu_p(t_1 + t_2) \geq \min\{\nu_p(t_1), \nu_p(t_2)\}$, multiplying two root elements associated to the same root does not decrease the valuation.
\end{remark}

The proof of the following lemma is almost identical to \cite[Lemma 3]{Tot05}; we reproduce it here for completeness. 
\begin{lemma}\label{isolating:bad:valuation}
Let $ab = \prod_{r\in\Phi_1} x_r(t_r)\, x_\hr(t_\hr) \in A(O')B(O')$.
\begin{enumerate}[label=(\alph*)]
  \item\label{isolating:height1} Let $\nu_p(a)=(n,s)<0$. Then there are a unit $\varepsilon\in O^\times$ and elements $t_s',t_\hr'\in O'$ with $\nu_p(t_s')>n$ and $\nu_p(t_\hr')\geq\min\{2n,\nu_p(t_\hr)\}$ such that
  \[
  ab = x_s(\varepsilon p^n)\, x_s(t_s') \prod_{r\in\Phi_1\setminus\{s\}} x_r(t_r)\; x_\hr(t_\hr').
  \]
  \item\label{isolating:height2} Let $\nu_p(t_\hr)=N<0$. Then there are a unit
  $\varepsilon\in O^\times$ and an element $t_\hr'\in O'$ with
  $\nu_p(t_\hr')>N$ such that
  \[
  ab = x_\hr(\varepsilon p^N)\, x_\hr(t_\hr') \prod_{r\in\Phi_1} x_r(t_r).
  \]
\end{enumerate}
\end{lemma}

\begin{proof}
We prove \ref{isolating:height1}. The proof of \ref{isolating:height2} is the same and simpler: since $x_\hr$ commutes with every other factor, bringing it to the first position creates no new terms.

Let $\nu_p(a)=(n,s)<0$. Write $t_{s} = k/p^m$ with $k \in O$. Since $n = \nu_p(t_{s}) < 0$, the valuation $\nu_p(k)$ is finite and $\nu_p(k) = m + n \geq 0$. Hence $k \in (p^{m+n})$, so $k = k' p^{m+n}$ with $k' \notin (p)$. Because $p$ is small, there is a unit $\varepsilon \in O^\times$ with $k' \equiv \varepsilon \pmod{p}$, say $k' = \varepsilon + p k''$ for some $k'' \in O$. Therefore
\[
t_{s} = k' p^{n} = \varepsilon p^{n} + k'' p^{n+1} = \varepsilon p^{n} + t_{s}',
\]
where $t_{s}' \coloneqq k'' p^{n+1}$ satisfies $\nu_p(t_{s}') \geq n+1 > n$.

We now bring $x_{s}(t_{s})$ to the first position and split it as
$x_{s}(t_{s}) = x_{s}(\varepsilon p^{n})\, x_{s}(t_{s}')$. Commuting $x_s(t_s)$
with another factor $x_r(t_r)$ of height $1$ produces, by
\Cref{swapping:elements}, at most one root element $x_\hr(u)$, with
$\nu_p(u) \geq \nu_p(t_s)+\nu_p(t_r) \geq 2n$, since $\nu_p(t_r)\geq n$ for every
$r \in \Phi_1$. These elements have height $2$, so they commute with all the
other factors and can be collected together with $x_\hr(t_\hr)$ into a single
element $x_\hr(t_\hr')$ with $\nu_p(t_\hr') \geq \min\{2n, \nu_p(t_\hr)\}$. 
\end{proof}

Here it is important that the prime $p$ is small, otherwise $\varepsilon$ might not be a unit. Since in our case it is a unit, we have $\nu_p(\varepsilon p^n)=\nu_p(p^n)=n$.

The next lemma gives information about the $\mathrm{SL}_2$-relations between proportional roots, hence it does not depend on the ambient group.

\begin{lemma}[{Tóth, \cite[Lemma 5]{Tot05}}] \label{proportional:roots}
Let $r$ be a positive root, $t_{-r} \in O$, and $t_r \in O'$ such that $\nu_p(t_{-r}) + \nu_p(t_r) >0$. Then
\[
x_{-r}(t_{-r}) x_r(t_r)= x_r(t_r') \gamma,
\]
for some $t'_r \in O'$ with $\nu_p(t_r') \geq \nu_p(t_r)$ and some $\gamma \in \G(O)$.
\end{lemma}

Throughout, $s$ denotes the smallest root at which the minimal valuation
$n$ of $a$ is attained. We will use repeatedly that this means
\begin{equation} \label{minimality}
r < s \;\Longrightarrow\; \nu_p(t_r) > n,
\qquad\qquad
\nu_p(t_r) = n \;\Longrightarrow\; r \geq s .
\end{equation}
 
Since $\Phi$ is reduced, $2s$ is not a root, so $\hr \neq 2s$. Given $ab \in A(O')B(O')$ with $\nu_p(a) = (n,s)$, $n<0$, we distinguish two situations:
\begin{enumerate}[label=\textbf{(\Roman*)}]
  \item\label{case:I} $\nu_p(b) \geq (2n,\hr)$, and either $\nu_p(b) > (2n,\hr)$ or $\hr > 2s$;
  \item\label{case:II} $\nu_p(b) \leq (2n,\hr)$, and either $\nu_p(b) < (2n,\hr)$ or $\hr < 2s$.
\end{enumerate}
In \ref{case:I} we increase $\nu_p(a)$, keeping
$\nu_p(b) \geq (2n,\hr)$. In \ref{case:II} we
increase $\nu_p(b)$, keeping $\nu_p(a) \geq (n,s)$. 

\subsubsection*{Case \ref{case:I}}
 
Apply \Cref{isolating:bad:valuation}. Since $\mathcal{R}$ is a right
$\Z[\G(K)]$-module, checking that $ab \in \Z[\G(O)] + \mathcal{R}$ is
equivalent to checking that
$(x_s(\varepsilon p^n) + \rho ) \prod_{r \neq s} x_r(t_r) \in \Z[\G(O)] + \mathcal{R}$
for some $\rho \in \mathcal{R}$. Choosing
$\rho = e - x_s(\varepsilon p^n) - x_{-s}(\varepsilon^{-1} p^{-n})$, the
latter is
\begin{equation} \label{second:summand}
a'b' - x_{-s}(\eta p^{-n})\, a'b',
\end{equation}
where $\eta = \varepsilon^{-1}$, $\nu_p(a') > \nu_p(a)$ and $\nu_p(b') \geq (2n,\hr)$. 

\begin{remark}
Bringing $x_s(t_s)$ to the first position might create a term of height $2$ and valuation equal to $(2n,\hr)$. If this is the case, then $\hr>2s$, so we are still under the hypothesis of case \ref{case:I}. Indeed, let $x_r(t_r)$ be the root element of height $1$ that does not commute with $x_s(t_s)$ and produce a term $x_{s+r}(t_{s+r}) = x_\hr(t_\hr)$ with $\nu_p(t_\hr)=2n$. Then $\nu_p(t_r)=n$, hence $r>s$ by \ref{minimality}. From this, $\hr=s+r > s+s = 2s$ by additivity of the chosen order in \Cref{root:order}.
\end{remark}

It is sufficient to check the claim for the second summand of \ref{second:summand}, which equals
\begin{align*}
x_{-s}(\eta p^{-n}) \prod_{\substack{r \in \Phi_{>0} \\ r \neq s}} x_r(t_r)\, x_s(t_s)
&= \left[ \prod_{\substack{r \in \Phi_{>0} \\ r \neq s}} x_{-s}(\eta p^{-n}) x_r(t_r) x_{-s}(-\eta p^{-n}) \right]
x_{-s}(\eta p^{-n})x_s(t_s) \\
&= \left[ \prod_{\substack{r \in \Phi_{>0} \\ r \neq s}} X(r;s)\right]
x_{-s}(\eta p^{-n})x_s(t_s),
\end{align*}
where by the Chevalley Commutator Formula
\[
X(r;s)=x_r(t_r) \prod_{\substack{i,j>0\\ ir-js\in\Phi}}
x_{ir-js}\left(C^{r,-s}_{i,j}\,t_r^{\,i}\,\eta^{\,j}\,p^{-jn}\right).
\]
To simplify the notation, we write
$X(r;s) = x_r(t_r) \prod_{q=ir-js} x_q(t_q)$. The last two factors give
$x_{-s}(\eta p^{-n})x_s(t_s) = x_s(t_s')\gamma$ for some
$\gamma \in \G(O)$ and $t'_s \in O'$ with $\nu_p(t'_s) >n$ by \Cref{proportional:roots}, since
$\nu_p(t_s) > n$ and hence $\nu_p(\eta p^{-n}) + \nu_p(t_s) > 0$.
 
We now split into cases depending on the height of $q = ir-js$ and the
valuation of $t_q$. If $r$ has height $1$, the possible types are the
following.
 
\begin{table}[H]
\centering
\begin{tabular}{c|c|c}
$(i,j)$ & $\operatorname{ht}(q)$ & $\nu_p(t_q)$\\
\hline
$(3,1)$ & $2$ & $\geq 2n$\\
$(3,2)$ & $1$ & $\geq n $\\
$(2,1)$ & $1$ & $\geq n$\\
$(1,1)$ & $0$ & $\geq 0$\\
$(1,2)$ & $-1$ & $\geq -n$\\
$(2,3)$ & $-1$ & $\geq -n$\\
$(1,3)$ & $-2$ & $\geq -2n$
\end{tabular}
\caption{Types for $r$ of height $1$.}
\label{tab:height1}
\end{table}
 
If instead $r = \hr$ has height $2$, the possible types are the following.
 
\begin{table}[H]
\centering
\begin{tabular}{c|c|c}
$(i,j)$ & $\operatorname{ht}(q)$ & $\nu_p(t_q)$\\
\hline
$(1,1)$ & $1$ & $\geq n$\\
$(2,3)$ & $1$ & $\geq n$\\
$(1,2)$ & $0$ & $\geq 0$\\
$(1,3)$ & $-1$ & $\geq -n$
\end{tabular}
\caption{Types for $r=\hr$ of height $2$.}
\label{tab:height2}
\end{table}
 
The bounds follow from $\nu_p(t_q) \geq i\,\nu_p(t_r) - jn$ together with
$\nu_p(t_r) \geq n$, respectively $\nu_p(t_\hr) \geq 2n$, and the list of
possible pairs $(i,j)$ follows from the possible rank $2$ systems spanned
by $r$ and $-s$. Note that if the root system is of type $G_2$, all the
cases of \Cref{tab:height1} are possible.
 However, if $\Phi \in \{F_4, E_8\}$, several combinations cannot appear.
For example, in \Cref{tab:height1}, the case $\operatorname{ht}(q)=-1$ can happen
only if $(i,j)=(1,2)$ and $\Phi = F_4$, and in this case $q$ is long.
This again follows from looking at the root system of rank $2$ spanned by
$r$ and $-s$: the combination $r-2s$ can be a root only if that subsystem
is of type $B_2$ and $s$ is short, and then $r-2s$ is long.
 
Observe that an equality in \Cref{tab:height1} forces $\nu_p(t_r)=n$,
hence $r > s$ by \eqref{minimality}, whereas an equality in
\Cref{tab:height2} forces $\nu_p(t_\hr)=2n$, hence $\hr > 2s$ by the
hypothesis of case \ref{case:I}. These two observations will be used in the proof of the following lemma.
 
\begin{lemma} \label{increasing:val:a}
Let $ab \in A(O')B(O')$ with $\nu_p(a) >(n,s)$, $\nu_p(b) \geq (2n,\hr)$, and either $\nu_p(b) > (2n,\hr)$ or $\hr>2s$. Let $x_q(t_q)$ be
as above.
Then there are elements $a'$, $b'$, $\gamma$ in $A(O')$, $B(O')$,
$\G(O)$ respectively, with $\nu_p(a')>(n,s)$ and
$b'$ satisfying the hypothesis of Case \ref{case:I}, such that
\[
x_q(t_q)\,ab=a'b'\gamma.
\]
\end{lemma}

To streamline the proof we show only once how to prove that $b'$ satisfies the hypothesis of Case \ref{case:I}; the other cases are treated analogously.
 
\begin{proof}
\noindent{\bfseries\boldmath Case $\operatorname{ht}(q)=2$.}
Since $q$ has height $2$, the root element $x_q(t_q)$ commutes with $a$ and $b$. We set $a'=a$ and $b'= x_q(t_q)b$, so that $\nu_p(b')\geq(2n,\hr)$. If $\nu_p(b')=(2n,\hr)$, then either $\nu_p(t_\hr)=2n$, and $\hr>2s$ by hypothesis, or $\nu_p(t_q)=2n$. In the latter case $q=3r-s$ with $\nu_p(t_r)=n$ by \Cref{tab:height1}, hence $r>s$ and $\hr=3r-s>2s$. In both cases $b'$ satisfies the hypothesis of case~\ref{case:I}.

\noindent{\bfseries\boldmath Case $\operatorname{ht}(q)=1$.}
We have $\nu_p(t_q) \geq n$, and if the inequality is strict, the claim follows. If $\nu_p(t_q) = n$, we must show
that $q > s$. If $q$ comes from \Cref{tab:height1}, then $r>s$ and
$q = 2r-s$ or $q = 3r-2s$, so $q > s$. If $q$ comes
from \Cref{tab:height2}, then $\hr>2s$ and $q = \hr-s$ or $q = 2\hr-3s$. In either case we have again $q >s$.
 
\noindent{\bfseries\boldmath Case $\operatorname{ht}(q)=0$.}
Using the Chevalley Commutator Formula, we have
\begin{align*}
x_q(t_q) \left[ \prod_{r \in \Phi_{>0}} x_r(t_r) \right] &=
 \left[ \prod_{r \in \Phi_{>0}} x_q(t_q) x_r(t_r) x_q(-t_q) \right] x_q(t_q) \\
&= \left[\prod_{r \in \Phi_{>0}} x_r(t_r) \prod_{iq+jr \in \Phi} x_{iq+jr}\!\left(C^{q,r}_{i,j} t_q^{\,i} t_r^{\,j}\right) \right] x_q(t_q).
\end{align*}

Since $\operatorname{ht}(q)=0$, a new factor $x_{iq+jr}$ occurs only for
$r\in\Phi_1$, and it has height $j\in\{1,2\}$.

Its valuation is at least $i\,\nu_p(t_q)+j\,\nu_p(t_r)\geq jn$, and equality forces $\nu_p(t_q)=0$ and $\nu_p(t_r)=n$, hence $q>0$ and $r>s$. Therefore, a new factor of height $1$ and valuation $n$ lies at $iq+r>s$, so $\nu_p(a')>(n,s)$.

\noindent{\bfseries\boldmath Case $\operatorname{ht}(q)=-1$, $\Phi=G_2$.}
We first record that
\begin{equation} \label{pairing}
\nu_p(t_q) + \nu_p(t_{-q}) > 0 .
\end{equation}
Indeed, if $\nu_p(t_q) > -n$ this is clear, since $\nu_p(t_{-q}) \geq n$.
If $\nu_p(t_q) = -n$, then $q > -s$, so $-q < s$ and $\nu_p(t_{-q}) > n$ by minimality of $s$, whence \eqref{pairing}.
 
Consider first the cases where $q$ is long, i.e.\
$q \in \{-3 \alpha - \beta, -\beta \}$. We proceed by induction on the
number of non-trivial root elements in $a$. If $a$ is the identity, then
\[
x_q(t_q)\,x_\hr(t_\hr) = x_q(t_q)\,x_\hr(t_\hr)\,x_q(t_q)^{-1}\,x_q(t_q)
= x_{q+\hr}\bigl(C^{q,\hr}_{1,1}t_q t_\hr\bigr)\,x_\hr(t_\hr)\,x_q(t_q),
\]
since $q+\hr$ is the only root of the form $iq+j\hr$ with $i,j>0$. The root
element $x_{q+\hr}(C^{q,\hr}_{1,1}t_qt_\hr)$ has height $1$ and valuation at least
$\nu_p(t_q)+\nu_p(t_\hr)\geq -n+2n=n$; if equality holds, then
$\nu_p(t_q)=-n$ and $\nu_p(t_\hr)=2n$, so $q>-s$ and $\hr>2s$, whence $q+\hr>s$. Moreover, $x_q(t_q) \in \G(O)$, so the claim follows. For the inductive step, consider a root $r_0 \neq -q$ in $ab$ such that $t_{r_0} \neq 0$. If such a choice does not
exist, then we are dealing with proportional root elements, and we may
apply \Cref{proportional:roots} by \eqref{pairing}. Otherwise, shifting
$x_{r_0}(t_{r_0})$ to the second position, we obtain
\begin{align*}
x_q(t_q) ab &= \left[ x_q(t_q)x_{r_0}(t_{r_0}) \right] \prod_{r_i \neq r_0} x_{r_i}(t_{r_i}) \\
&= \left[ x_q(t_q) x_{r_0}(t_{r_0}) x_q(-t_q) \right]  x_q(t_q) \prod_{r_i \neq r_0} x_{r_i}(t_{r_i}) \\
&= \left[ x_q(t_q) x_{r_0}(t_{r_0}) x_q(-t_q) \right] a' b' \gamma \\
&= \left[ x_{r_0}(t_{r_0}) \prod_{iq+jr_0 \in \Phi} x_{iq + jr_0}\!\left(C^{q,r_0}_{i,j} t_q^{\,i} t_{r_0}^{\,j}\right) \right] a' b' \gamma,
\end{align*}
where the last two equalities follow from the inductive hypothesis and
the Chevalley Commutator Formula respectively. The combinations
$iq + jr_0$ can only be roots of height greater than or equal to $0$,
and these are dealt with similarly as above.

For short $q$ an explicit computation is needed, because the inductive
argument used for long roots does not terminate: the commutator of
$x_{-\alpha-\beta}$ and $x_{\beta}$ produces a factor
$x_{-2\alpha-\beta}$, while the commutator of $x_{-2\alpha-\beta}$ and
$x_{3\alpha+\beta}$ produces a factor $x_{-\alpha-\beta}$, so the induction would keep alternating between the two short roots of height $-1$. 
 
Suppose $q = -\alpha -\beta$. If $ab$ contains a root element which is
associated neither to $\beta$ nor to $-q$, we repeat the inductive
procedure carried out in the previous step: we already handled the
negative terms that might appear. If instead we are in the case
$x_{-\alpha - \beta}(t_{-\alpha - \beta}) x_{\beta}(t_\beta) x_{\alpha + \beta}(t_{\alpha + \beta})$,
an explicit computation shows the result. In the following,
$\gamma^{(i)}$ indicate elements of $\G(O)$ for every $i=1,\dots, 6$. We
do not write the explicit coefficients to lighten the notation; every root
element $x_{q'}(t)$ appearing below satisfies $\nu_p(t)\geq n\operatorname{ht}(q')$.
\begin{align*}
&\x{-\alpha-\beta}\x{\beta}\x{\alpha+\beta}  \\
={}&\x{-\alpha-\beta}\x{\beta} \x{-\alpha-\beta}^{-1}
   \x{-\alpha-\beta}\x{\alpha+\beta}  \\
={}&\x{\beta}\x{-\alpha}\x{-2\alpha-\beta}\x{-3\alpha-2\beta}
   \x{-3\alpha-\beta} \x{\alpha+\beta}\, \gamma^{(1)} \\
={}&\x{\beta}\x{-\alpha}\x{-2\alpha-\beta}\x{-3\alpha-2\beta}
   \x{\alpha+\beta}\, \gamma^{(2)}  \\
={}&\x{\beta}\x{-\alpha}\x{-2\alpha-\beta}\x{\alpha+\beta}\x{-2\alpha-\beta} \x{-\alpha}\x{\beta}\x{-3\alpha-2\beta}\, \gamma^{(3)}  \\
={}&\x{\beta}\x{-\alpha}\x{-2\alpha-\beta}\x{\alpha+\beta}\x{-2\alpha-\beta} \x{\beta}\, \gamma^{(4)}  \\
={}&\x{\beta}\x{-\alpha}\x{-2\alpha-\beta}\x{\alpha+\beta}
   \x{-2\alpha-\beta}^{-1} \x{-2\alpha-\beta}\x{\beta}  \, \gamma^{(5)}  \\
={}&\x{\beta}\x{-\alpha}\x{\alpha+\beta}\x{-\alpha}\x{\beta}
   \x{-3\alpha-\beta}\x{\beta}\, \gamma^{(6)}.
\end{align*}
The terms appearing in the last row are treated by as cases above. Note that we are still under the hypothesis of the lemma: the valuation of the height $1$ terms is strictly bigger than $(n,s)$.
 
The case $q = -2\alpha-\beta$ is handled as the case $q = -\alpha-\beta$.
This time the root that needs further checks is $3\alpha+\beta$.

\noindent{\bfseries\boldmath Case $\operatorname{ht}(q)=-1$, $\Phi \in \{F_4, E_8\}$.}
We proceed by induction as we did for $G_2$, and the treatment of
proportional root elements is the same, since \eqref{pairing} holds here
as well. In $F_4$ a root $q$ with height $-1$ arising in this case is necessarily long, and in
the inductive step there are no longer any factors of negative height
after applying the Chevalley Commutator Formula. For $E_8$ no factor of
negative height arises either: the only factor produced is $x_{q+r_0}$, of height $\operatorname{ht}(r_0) - 1 \geq 0$. These facts
follow by looking at which root systems of rank $2$ embed into $\Phi$:
for $E_8$ this is the case only for $A_1 \times A_1$ and $A_2$, since it
is simply laced.
 
\noindent{\bfseries\boldmath Case $\operatorname{ht}(q)=-2$.}
Note that this can happen only when $\Phi = G_2$. We do the same inductive
procedure as we did for $\operatorname{ht}(q)=-1$. For the base step, we claim that
\[
\nu_p(t_{-\hr}) + \nu_p(t_\hr) > 0,
\]
so that \Cref{proportional:roots} applies. Indeed, if $\nu_p(t_{-\hr}) > -2n$ this is clear, since $\nu_p(t_\hr) \geq 2n$. If instead $\nu_p(t_{-\hr}) = -2n$, then by \Cref{tab:height1} we have $q = r-3s$ for some $r \in \Phi_1$ with $\nu_p(t_r) = n$, hence $r>s$ and
\[
-\hr = q = r-3s > s-3s = -2s,
\]
that is $\hr < 2s$. By the hypothesis of case \ref{case:I} this forces
$\nu_p(b) > (2n,\hr)$, so $\nu_p(t_\hr) > 2n$ and the claim follows.

 In the inductive step we reduce to studying the possible types of $iq+jr_0$, which have height either $-1$, $0$ or $1$,
and all these cases have already been dealt with.

\end{proof}

\subsubsection*{Case \ref{case:II}}
 
Here $\nu_p(b) = (N,\hr) \leq (2n,\hr)$, and either the inequality is strict or  $\nu_p(b) = (2n,\hr)$ and $\hr<2s$. The computation above cannot be
carried out: for instance, if $\nu_p(b)=(2n,\hr)$, for $(i,j) = (1,1)$ in \Cref{tab:height2} one
gets
\[
q = \hr-s < 2s-s = s,
\]
and a factor of height $1$ and valuation less than $(n,s)$ would be created.
We therefore isolate the bad valuation using $x_\hr(t_\hr)$ instead of
$x_s(t_s)$, increasing $\nu_p(b)$ and leaving $\nu_p(a)$
unchanged.
 
Namely, we bring $x_\hr(t_\hr)$ to the first position, apply
\Cref{isolating:bad:valuation} and add a suitable element of
$\mathcal{R}$, reducing to proving the statement for
\begin{align*}
x_{-\hr}(\eta p^{-N}) \left[ \prod_{r \in \Phi_1} x_r(t_r) \right] x_\hr(t_\hr)
&= \left[ \prod_{r \in \Phi_1} X(r;\hr)\right]
x_{-\hr}(\eta p^{-N})x_\hr(t_\hr'),
\end{align*}
where $\nu_p(t_r)\geq n$, $\nu_p(t_\hr') > N$ and
\[
X(r;\hr)=x_r(t_r) \prod_{\substack{i,j>0\\ ir-j\hr \in \Phi}}
x_{ir-j\hr}\left(C^{r,-\hr}_{i,j}\,t_r^{\,i}\,\eta^{\,j}\,p^{-jN}\right).
\]
The last two factors give $x_\hr(t_\hr')\gamma$ for some $\gamma \in \G(O)$
by \Cref{proportional:roots}, since
$\nu_p(\eta p^{-N}) + \nu_p(t_\hr') > -N + N = 0$. Writing as before
$X(r;\hr)=x_r(t_r) \prod_{q=ir-j\hr} x_q(t_q)$, the possible types of $x_q(t_q)$ are
the following.
 
\begin{table}[H]
\centering
\begin{tabular}{c|c|c}
$(i,j)$ & $\operatorname{ht}(q)$ & $\nu_p(t_q)$\\
\hline
$(3,1)$ & $1$ & $\geq n$\\
$(2,1)$ & $0$ & $\geq 0$\\
$(1,1)$ & $-1$ & $\geq -n$\\
$(3,2)$ & $-1$ & $\geq -n $
\end{tabular}
\caption{Types for $q=ir-j\hr$.}
\label{tab:height2bis}
\end{table}
 
\begin{lemma} \label{increasing:val:b}
Let $N \leq 2n$, and assume $\hr<2s$ if $N=2n$. Let $ab \in A(O')B(O')$ with
$\nu_p(a) \geq (n,s)$ and $\nu_p(b) > (N,\hr)$, and let $x_q(t_q)$ be as above. Then there are elements $a'$, $b'$, $\gamma$ in $A(O')$, $B(O')$, $\G(O)$ respectively, with
$\nu_p(a')\geq(n,s)$ and $\nu_p(b') > (N,\hr)$, such that
\[
x_q(t_q)\,ab=a'b'\gamma.
\]
\end{lemma}
 
\begin{proof}
The proof is analogous to \Cref{increasing:val:a}. For instance, for $\operatorname{ht}(q)=1$ the claim follows by absorbing $x_q(t_q)$ into $a$; for $\operatorname{ht}(q)=0$ the terms arising by the Chevalley Commutator Formula have height $1$ or $2$ and valuation greater than or equal to $(n,s)$ and strictly greater than $(N,\hr)$ respectively. Those of height $2$ are root elements $x_\hr(t_\hr)$ with $\hr=iq+2r$ and $\nu_p(t_\hr)\geq i\nu_p(t_q)+2\nu_p(t_r)\geq 2n$. Equality would force $\nu_p(t_r)=n$, hence $r\geq s$ and $\hr=iq+2r>2s$, contradicting $\hr<2s$. Therefore $\nu_p(t_\hr)>2n$, and $\nu_p(b)>(N,\hr)$.

The only combination not yet discussed occurs in $F_4$ for $\operatorname{ht}(q)=-1$, where $q = r-\hr$ is short whenever $r$ is short; conjugating by such a $q$, the only factors of negative height are those of the form $2q+r_0$ with $r_0$ long of height $1$, and these are
long roots of height $-1$, which have already been treated.
\end{proof}
 
 \begin{remark}
The following observation clarifies what happens when $\nu_p(b)=(2n,\hr)$ and the importance of the assumption $\hr<2s$. We start with an element $ab \in A(O')B(O')$ with $\nu_p(a) = (n,s)$, $\nu_p(b) \leq (2n,\hr)$, we increase the valuation of $b$ until it reaches $(2n,\hr)$. If $\hr>2s$, we apply case \ref{case:I} and increase the valuation of $a$. If $\hr<2s$, we increase the valuation of $b$ again and then apply case \ref{case:I}. To do so, we need to show that swapping two root elements of height $1$ and valuation greater than or equal to $(n,s)$ does not produce a term of height $2$ and valuation equal to $(2n,\hr)$. Let $a'b' \in A(O')B(O')$ the element obtained from $ab$ after iterating case \ref{case:II}. Let $x_{r_1}(t_{r_1})$ and $x_{r_2}(t_{r_2})$ root elements appearing as factors in $a'$ such that $\nu_p(t_{r_1})=\nu_p(t_{r_2}) = n$. The latter implies that $r_i \geq s$ for $i=1,2$, and the inequality is strict for at least one of the two roots. Thus, $r_1 + r_2 > 2s$, and since we assumed $\hr< 2s$, the two root elements commute. 
\end{remark}

Finally, suppose that all factors of $a$ have nonnegative valuation, while
$\nu_p(b)=(N,\hr)$ with $N<0$. Then we proceed as in Case~\ref{case:II} with
$n=0$, and the argument goes through verbatim.
 
\subsubsection*{Conclusion}
 
Let us summarize what we did in this section. We want to prove
\Cref{local:thm}, which says that an element of $U(O')$ can be written as
something in $\Z[\G(O)] + \mathcal{R}$. Thanks to \Cref{height 0}, we can
instead work with an element $ab \in A(O')B(O')$.

We have written such an element as $ab \equiv a'b' + a''b''\gamma \pmod{\mathcal{R}}$, and according to which of the two situations \ref{case:I}, \ref{case:II} occurs, we strictly increased $\nu_p(a)$ or $\nu_p(b)$. Iterating this process, we eventually obtain terms lying in $\Z[\G(O)] + \mathcal{R}$.

\subsection{The global statement} \label{global:section}
In this section we prove \Cref{integral:generation}, i.e. that the Steinberg module is generated by integral apartment classes, as done in \cite[Section 6]{Tot05}.

Recall the setup. Let $\G$ be a semisimple Chevalley--Demazure group scheme of type either $G_2$, $F_4$ or $E_8$, and let $K$ be field of fractions of a Euclidean domain $O$. The Steinberg module of $\G(K)$ is denoted by $\mathrm{St}$, and its fundamental apartment class by $[\Sigma]$.

By the Solomon-Tits~\Cref{solomon:tits:thm}, we have the surjection
\[
\pi : \Z[\G(K)] \twoheadrightarrow \mathrm{St},\qquad g \mapsto [g]=[\Sigma]\cdot g,
\]
and $\mathcal{R} \subseteq \ker(\pi)$ (see \Cref{R:ker}). If we prove that $\Z[\G(K)] = \Z[\G(O)] + \mathcal{R}$, we obtain
\[
\mathrm{St}=\pi\big(\Z[\G(K)]\big)=\pi\big(\Z[\G(O)]\big)=[\Sigma]\cdot\Z[\G(O)],
\]
i.e. the Steinberg module is generated by integral apartment classes. In the rest of this section we prove this equality.

Let $p \in O$ be a prime. For $t \in O$, recall that we defined
\[
\nu_p(t)\coloneqq\sup\{n \in \N : p^n\mid t\}
\]
with the convention $\nu_p(0)=\infty$. For $t/s \in K$, we set $\nu_p(t/s) \coloneqq \nu_p(t) - \nu_p(s)$.

For a prime $p \in O$, let $O_p = \{t \in K \mid \nu_p(t) \geq 0\}$ denote the ring of \emph{$p$-integral} elements of $K$. We say $g \in \G(K)$ is \emph{$p$-integral} if $g \in \G(O_p)$, and set
\[
  \mathfrak{S}_g = \{\, p \text{ prime} \mid g \notin \G(O_p) \,\}.
\]
It is a finite set, and $\mathfrak{S}_g = \varnothing$ if and only if $g \in \G(O)$. Indeed, from the ring equality $O=\bigcap_p O_p$ we have $\G(O) = \bigcap_p \G(O_p)$. 

Let $d$ be a Euclidean function on $O$ satisfying $d(x) \geq d(xy)$ for every $x,y \in O$ and fix $g \in \G(K)$. Let $n=\max\{d(p) : p\in \mathfrak{S}_g\}$ with the convention that $n=-\infty$ if $\mathfrak{S}_g=\emptyset$, and let $m$ be the number of primes $p \in \mathfrak{S}_g$ with $d(p)=n$. We endow the pairs $(n,m)$ with the lexicographic order and we prove that $g \in \Z[\G(O)] + \mathcal{R}$ by induction. The base step corresponds to the case where $g$ is already integral.

For the inductive step, it will be useful to reduce to working with elements in $U$ so that we can apply \Cref{local:thm}. To do so, we recall the \emph{Iwasawa decomposition} for a Chevalley group.
\begin{lemma}[{\cite[Theorem 18, Chapter 8]{Ste16}}] \label{iwasawa:decomposition}
Let $O$ be a Euclidean domain with fraction field $K$ and let $\G$ be a Chevalley--Demazure group scheme with fixed maximal torus $T$ and unipotent radical $U$ of the Borel. Then, every $g \in \G(K)$ can be decomposed as 
\[
g=hu \gamma \in T(K)U(K) \G(O).
\]
Moreover, for every such decomposition and every prime $p$ of $O$, if $g\in\G(O_p)$, then $u\in U(O_p)$.
\end{lemma}

Let $g = hu\gamma$ with $h \in T(K)$, $u \in U(K)$ and $\gamma \in \G(O)$. If we prove $u \gamma \in \Z[\G(O)] + \mathcal{R}$, then the claim follows. Indeed:  
\[
u \gamma - (e-h) \cdot u\gamma = hu \gamma = g \in \Z[\G(O)] + \mathcal{R}.
\]
Since $\gamma$ is already integral, we only need to show the claim for $u \in U(K)$. Let $\{p_1, \dots, p_m\} = \{p \in \mathfrak{S}_g : d(p) = n \}$. We want to apply \Cref{local:thm}, so we need to work in a ring where the prime we are considering is small.

Let $O_{<n} \coloneqq O[c^{-1} : d(c) <n]$. Then $p_m$ is small in $O_{<n}$. For let $x/y \in O_{<n}$ with $x \in O$ and $y$ a unit of $O_{<n}$ such that $p_m \nmid x/y$. It suffices to check that $x$ can be represented by a unit. Since $p_m \nmid x$, we can write $x = q p_m + r$ with $r \neq 0$ and $d(r) < d(p_m)=n$. Then $x \equiv r \mod{p_m}$ and $r$ is a unit of $O_{<n}$.

Let $\tilde{O} \coloneqq  O_{<n}[p_1^{-1}, \dots , p_{m-1}^{-1}]$. Then an analogous argument to the one above shows that $p_m$ is still small in $\tilde{O}$. Moreover, $u \in \Z[U(\tilde{O}[p_m^{-1}])]$ since $u$ is $p$-integral where $g$ is. By \Cref{local:thm}, the element $u$ can be written as $u=\sum_i c_i \gamma_i + \rho$ with $c_i \in \Z$, $\gamma_i \in \G(\tilde{O})$ and $\rho \in \mathcal{R}$. 
For every $\gamma_i$, we have $\mathfrak{S}_{\gamma_i}\subseteq\{p:d(p)<n\}\cup\{p_1,\dots,p_{m-1}\}$. Applying the inductive hypothesis to each of the $\gamma_i$'s, we obtain the claim.

  \bibliographystyle{amsalpha}
  \bibliography{references}
\end{document}